\documentclass{llncs}
\usepackage{amsmath}
\usepackage{amssymb}
\usepackage{xspace}
\usepackage{color}

\newenvironment{Proof}{\begin{proof}}{\qed\end{proof}}

\newcommand{\field}{\mathbb{F}}
\newcommand{\sub}[1]{[ #1 ]}
\newcommand{\sdif}{\mathop{\mathrm{\Delta}}}
\newcommand{\dmatroid}{$\Delta$-matroid\xspace}
\newcommand{\dmatroids}{$\Delta$-matroids\xspace}
\newcommand{\dual}{\mathop{\bar{*}}}

\newcommand{\vfsafe}{vf-safe\xspace}

\renewcommand{\emptyset}{\varnothing}

\title{Quaternary matroids are vf-safe}
\author{Robert Brijder\inst{1} \and Hendrik Jan Hoogeboom\inst{2}}
\institute{Hasselt University and Transnational University of Limburg, Belgium \email{robert.brijder@uhasselt.be} \and Leiden Institute of Advanced Computer Science,\\
Leiden University, The Netherlands\\ \email{hoogeboom@liacs.nl}}

\begin{document}

\maketitle

\begin{abstract}
Binary delta-matroids are closed under vertex flips, which consist of the natural operations of twist and loop complementation. In this note we provide an extension of this result from $GF(2)$ to $GF(4)$. As a consequence, quaternary matroids are ``safe'' under vertex flips (\vfsafe for short). As an application, we find that the matroid of a bicycle space of a quaternary matroid is independent of the chosen representation. This extends a result of Vertigan
[J.\ Comb.\ Theory B (1998)] concerning the bicycle dimension of quaternary matroids.
\end{abstract}


\section{Preliminaries}

\subsection{Notation and terminology}
For finite sets $U$ and $V$, a $U \times V$-matrix $A$ (over
some field $\field$) is a matrix where the rows are indexed by
$U$ and the columns by $V$, i.e., $A$ is formally a function $U
\times V \rightarrow \field$. Hence, the order of the
rows/columns is not fixed (i.e., interchanging rows or columns
is mute). For $X \subseteq U$ and $Y \subseteq V$, the
submatrix of $A$ induced by $X$ and $Y$ is denoted by $A[X,Y]$. We often abbreviate $A[X,X]$ by $A[X]$. Let $A$ be a $V \times V$-matrix and $I_X$ for $X \subseteq V$ be the $V \times V$-matrix
where the the diagonal entries corresponding to $X$ are $1$ and all other entries are $0$. We abbreviate $A+I_X$ by $A+X$.

\subsection{Principal pivot transform}
Let $\alpha$ be an automorphism of a field $\field$. By abuse of notation, we extend $\alpha$ point-wise to vectors, matrices, and subspaces over $\field$. Hence for a $V \times V$-matrix $A = (a_{i,j})_{i,j \in V}$, we let $\alpha(A) = (\alpha(a_{i,j}))_{i,j \in V}$. Moreover, for subspace $L \subseteq \field^V$, we let $\alpha(L) = \{ \alpha(v) \mid v \in L\}$.

Let $A$ be a $V \times V$-matrix (over an arbitrary field), and
let $X \subseteq V$ be such that $A\sub{X}$ is nonsingular,
i.e., $\det A\sub{X} \neq 0$. The \emph{principal pivot
transform} (or \emph{PPT} for short) of $A$ on $X$, denoted
by $A*X$, is defined as follows, see \cite{tucker1960}. If $A
= \bordermatrix{
& X & V\setminus X \cr
X & P & Q \cr
V\setminus X & R & S
}$, then
$$
A*X = \bordermatrix{
& X & V\setminus X \cr
X & P^{-1} & -P^{-1} Q \cr
V\setminus X & R P^{-1} & S - R P^{-1} Q
}.
$$
Matrix $(A*X)\setminus X = S - R P^{-1} Q$ is called the
\emph{Schur complement} of $X$ in $A$. Hence, $A*X$ is defined
iff $A\sub{X}$ is nonsingular.

It is easy to verify (by the above definition of PPT) that $-(A*X)^T = (-A^T)*X$ for all $X \subseteq V$ with $A[X]$ nonsingular. As a consequence, if $A$ is skew-symmetric, then $A*X$ is skew-symmetric.

\begin{lemma} \label{lem:automorph_PPT}
Let $A$ be a $V \times V$-matrix over some field $\field$, and let $\alpha$ be an automorphism of $\field$. If $X \subseteq V$ is such that $A[X]$ is nonsingular, then $\alpha(A*X) = \alpha(A)*X$.
\end{lemma}
\begin{Proof}
Obviously, for a nonsingular matrix $P$, $P P^{-1} = I$, where $I$ is the identity matrix (of suitable size). Therefore, $I = \alpha(P P^{-1}) = \alpha(P) \alpha(P^{-1})$, and so $\alpha(P^{-1}) = \alpha(P)^{-1}$.
If $A
= \bordermatrix{
& X & V\setminus X \cr
X & P & Q \cr
V\setminus X & R & S
}$, then in both cases we obtain:
$$
\bordermatrix{
& X & V\setminus X \cr
X & \alpha(P)^{-1} & -\alpha(P)^{-1} \alpha(Q) \cr
V\setminus X & \alpha(R) \alpha(P)^{-1} & \alpha(S) - \alpha(R) \alpha(P)^{-1} \alpha(Q)
}.
$$
\end{Proof}

\begin{proposition}[\cite{tucker1960}]\label{prop:tucker}
Let $A$ be a $V \times V$-matrix, and let $X\subseteq V$ be
such that $A[X]$ is nonsingular. Then, for all $Y \subseteq
V$, $\det((A*X)[Y]) = \det(A[X \sdif Y]) / \allowbreak \det(A[X])$. In particular, $(A*X)[Y]$ is nonsingular iff $A[X \sdif Y]$ is nonsingular.
\end{proposition}


\subsection{Pivot and loop complementation on set systems}

A \emph{set system} (over $V$) is a tuple $M = (V,E)$ with $V$ a
finite set called the \emph{ground set} and $E \subseteq 2^V$ a
family of subsets of $V$. Set system $M$ is called \emph{proper} if $E \not= \emptyset$. We write simply $Y \in M$ to denote $Y \in E$. Let $\max(E)$ be the family of maximal sets in $E$ with respect to set inclusion, and let $\max(M) = (V,\max(E))$ be the corresponding set system.

We define, for $X \subseteq
V$, \emph{pivot} (also called \emph{twist}) of $M$ on $X$, denoted by $M * X$, as $(V,E *
X)$, where $E * X = \{Y \sdif X \mid Y \in E\}$. In case $X =
\{u\}$ is a singleton, we also write simply $M * u$. Moreover,
we define, for $u \in V$, \emph{loop complementation} of $M$ on
$u$, denoted by $M + u$, as
$(V,E')$, where $E' = E \sdif \{X \cup \{u\} \mid X \in E, u
\not\in X\}$. We assume left associativity of set system
operations. Therefore, e.g., $M+u *v$ denotes
$(M+u)*v$. It has been shown in \cite{BH/PivotLoopCompl/09} that pivot ${}*u$ and
loop complementation ${}+u$ on a common element $u \in V$ are
involutions (i.e., of order $2$) that generate a group $F_u$ isomorphic to
$S_3$, the group of permutations on $3$ elements. In
particular, we have ${}+u*u+u = {}*u+u*u$, which is the third
involution (in addition to pivot and loop complementation), and
is called the \emph{dual pivot}, denoted by $\dual $. The
elements of $F_u$ are called \emph{vertex flips}. We have, e.g., ${}+u*u = {}\dual
u+u = {}*u \dual  u$ and ${}*u+u = {}+u \dual u = {}\dual  u
*u$ for $u \in V$ (which are the two vertex flips in $F$ of
order $3$).

While on a single element the vertex flips behave as the group
$S_3$, they commute when applied on different elements. Hence, e.g., $M*u+v= M+v*u$ and $M\dual u+v= M+v\dual u$ when $u \not= v$. Also, $M+u+v= M+v+u$
and thus we (may) write, for $X = \{u_1,u_2,\ldots,u_n\} \subseteq
V$, $M+X$ to denote $M + u_1 \ldots + u_n$ (as the result is
independent on the order in which the operations ${}+u_i$ are
applied). Similarly, we define $M\dual X$ for $X \subseteq V$.

One may explicitly define the sets in $M*V$, $M+V$, and $M\dual
V$ as follows: $X \in M*V$ iff $V-X \in M$, and $X \in M+V$ iff
$|\{ Z\in M \mid Z\subseteq X \}|$ is odd. Dually, $X \in
M\dual V$ iff $|\{ Z\in M \mid X\subseteq Z \}|$ is odd. In
particular $\emptyset\in M\dual V$ iff the number of sets in
$M$ is odd.

Finally, it is observed in \cite{BH/PivotLoopCompl/09} that
$\max(M) = \max(M\dual X)$ for all $X\subseteq V$.

We will often use the results of this subsection without explicit
mention. Also we often simply denote the ground set of the set system under consideration by $V$.

\subsection{\vfsafe \dmatroids}

We assume that the reader is familiar with the basic notions of matroids, see, e.g., \cite{Oxley/MatroidBook-2nd}.

A \emph{\dmatroid} is a proper set system $M$ that satisfies
the \emph{symmetric exchange axiom}: For all $X,Y \in M$ and
all $u \in X \sdif Y$, either $X \sdif \{u\} \in M$ or there is
a $v \in X \sdif Y$ with $v \not= u$ such that $X \sdif \{u,v\}
\in M$ \cite{bouchet1987}. Note that \dmatroids are closed
under pivot, i.e., $M*X$ for $X \subseteq V$ is a \dmatroid
when $M$ is a \dmatroid. If we assume a
matroid $M$ is described by its basis, i.e., $M$ is the set system $(V,B)$ where $B$ is the
set of bases of $M$, then it is shown in
\cite[Proposition~3]{DBLP:conf/ipco/Bouchet95} that a matroid
$M$ is precisely a equicardinal \dmatroid. Hence \dmatroids form a generalization of matroids.

We say that a \dmatroid $M$ is \emph{vertex-flip-safe} (or \emph{\vfsafe} for short) if for any
sequence $\varphi$ of vertex flips (equivalently, pivots and
loop complementations) over $V$ we have that $M\varphi$ is a
\dmatroid. The family of \vfsafe \dmatroids is minor closed \cite{BH/NullityLoopcGraphsDM/10}.
We say that a family of \dmatroids is \emph{vf-closed} if the
family is closed under invertible vertex flips.
There are (delta-)matroids that are not \vfsafe, such as the $6$-point line $U_{2,6}$, $P_6$, and the non-Fano $F_7^-$. In fact, they are excluded minors for the family of \vfsafe \dmatroids.

Let $v \in \field^V$ be a vector. The support of $v$ is the set $X \subseteq V$ such that the entries of $X$ in $v$ are nonzero and entries of $V\setminus X$ in $v$ are zero. Let $L \subseteq \field^V$ be a subspace of $\field^V$. We denote by $M(L)$ the matroid with ground set $V$ such that for all $X \subseteq V$, $X$ is a circuit of $M(L)$ iff there is a $v \in L$ with support $X$ and $X$ is minimal with this property among the non-empty subsets of $V$. For a $X \times V$-matrix $A$, we denote by $M(A) = M(\ker(A))$, the matroid corresponding to the nullspace of $A$.

Let $A$ be a $V \times V$-matrix. We denote by $\mathcal{M}_A$ the set system $(V,D)$ where $D = \{ X \subseteq V \mid A[X] \mbox{ is nonsingular}\}$. It is shown in \cite{bouchet1987} (cf. Lemma~\ref{lem:alpha_symm_dmatroid}) that $\mathcal{M}_{A}$ is a \dmatroid if $A$ is skew-symmetric (i.e., $-A^T = A$). A \dmatroid $M$ is said to be \emph{representable} over $\field$, if $M = \mathcal{M}_A*X$ for some skew-symmetric $V \times V$-matrix $A$ and some $X \subseteq V$. It turns out that a matroid $M$ is representable in the \dmatroid sense iff $M$ is representable in the usual (matroid) sense. Moreover, if $A$ is skew-symmetric, then $\max(\mathcal{M}_{A})$ is a matroid represented by its bases and equal to $M(A)$ (this follows from the strong principal minor theorem \cite{Kodiyalam_Lam_Swan_2008}).

A \dmatroid is said to be \emph{binary} if it is representable over $GF(2)$. It is shown in \cite{BH/NullityLoopcGraphsDM/10} that the family of binary \dmatroids is vf-closed. Consequently, the class of \vfsafe \dmatroids contains the class of binary \dmatroids (and therefore also the class of binary matroids).

\section{$\alpha$-symmetry and delta-matroids}

Let $A$ be a $V \times V$-matrix over some field $\field$, and let $\alpha$ be an automorphism of $\field$. Then $A$ is called \emph{$\alpha$-symmetric} if $\alpha(-A^T) = A$. Note that if $A$ is $\alpha$-symmetric, then $\alpha(\alpha(x)) = x$ for all elements $x$ of $A$. Thus $\alpha$ behaves as an involution on the elements of $A$. As a consequence, if $A$ is $\alpha$-symmetric, then $A^T$ is $\alpha$-symmetric. Also note that $A$ is $\mathrm{id}$-symmetric with $\mathrm{id}$ the identity automorphism iff $A$ is skew-symmetric.


If $A$ is $\alpha$-symmetric and $X \subseteq V$ is such that $A[X]$ is nonsingular, then $A*X$ is $\alpha$-symmetric. Indeed, $\alpha(-(A*X)^T) = \alpha((-A^T)*X) = \alpha(-A^T)*X = A*X$, where in the second equality we use Lemma~\ref{lem:automorph_PPT}.

The next result is a straightforward extension of a result of \cite{bouchet1987} (the original formulation restricts to the case $\alpha = \mathrm{id}$).
\begin{lemma} [\cite{bouchet1987}]  \label{lem:alpha_symm_dmatroid}
Let $\alpha$ be an automorphism of some field $\field$, and let $A$ be a $\alpha$-symmetric $V \times V$-matrix over $\field$. Then $\mathcal{M}_{A}$ is a \dmatroid.
\end{lemma}
\begin{Proof}
Let $X,Y \in \mathcal{M}_{A}$ and $x \in X \sdif Y$. If entry $A*X[\{x\}]$ is nonzero, then by Proposition~\ref{prop:tucker}, $X \sdif \{x\} \in \mathcal{M}_{A}$ and we are done. Thus assume that $A*X[\{x\}]$ is zero. Since $A[Y]$ is nonsingular, $A*X[X\sdif Y]$ is nonsingular. Hence there is a $y \in X \sdif Y$ with entry $A*X[\{x\},\{y\}]$ nonzero (note that $x \neq y$). Since $A*X$ is $\alpha$-symmetric, $A*X[\{x,y\}]$ is of the form
$$
\bordermatrix{&x&y\cr
x & 0 & t_1 \cr
y & \alpha(-t_1) & t_2
}
$$
for some $t_1 \in \field \setminus \{0\}$ and $t_2 \in \field$. Thus $A*X[\{x,y\}]$ is nonsingular and $X \sdif \{x,y\} \in \mathcal{M}_{A}$.
\end{Proof}

We say that a \dmatroid $M$ is \emph{$\alpha$-representable} over $\field$, if $M = \mathcal{M}_A*X$ for some $\alpha$-symmetric $V \times V$-matrix $A$ and $X \subseteq V$. Note that this is a natural extension of the notion of representable from \cite{bouchet1987} which coincides with $\mathrm{id}$-representable.

A $V \times V$-matrix $A$ over $\field$ is called \emph{principally  unimodular} (PU, for short) if for all $Y \subseteq V$, $\det(A[Y]) \in \{0,1,-1\}$. Note that any $V \times V$-matrix over $GF(2)$ or $GF(3)$ is principally unimodular.

We now consider the field $GF(4)$. Let us denote the unique nontrivial automorphism of $GF(4)$ by $\mathrm{inv}$. Note that $\mathrm{inv}(x) = x^{-1}$ for all $x \in GF(4) \setminus \{0\}$, and thus $\mathrm{inv}$ is an involution.

\begin{theorem} \label{thm:inv_symm_PU}
Let $A$ be a $\mathrm{inv}$-symmetric $V \times V$-matrix over $GF(4)$. Then $A$ is a principally unimodular.
\end{theorem}
\begin{Proof}
Recall that $1 = -1$ in $GF(4)$. We have $\det(A) = \det(\mathrm{inv}(-A^T)) = \mathrm{inv}(\det(-A^T)) = \mathrm{inv}(\det(A))$. Thus $\det(A) \in \{0,1\}$.
\end{Proof}

\begin{remark}
The proof of Theorem~\ref{thm:inv_symm_PU} essentially uses that the field $\field$ under consideration is of characteristic $2$, i.e., $\field = GF(2^k)$ for some $k \geq 1$, \emph{and} $\field$ has an automorphism $\alpha$ such that $\alpha$ is an involution and $\alpha$ has only trivial fixed points (the set of fixed points form $GF(2)$). The automorphisms $\alpha$ of $GF(2^k)$ are of the form $x \mapsto x^{p^l}$, with $1 \le l \le k$, and $\alpha$ is an involution when either $k=1$ (and thus $l=1$) or both $k$ is even and $l = k/2$. Moreover, for $l = k/2$ and $k$ even, the corresponding automorphism $\alpha$ has only trivial fixed points iff $l=1$. Consequently, the proof of Theorem~\ref{thm:inv_symm_PU} only works for $\alpha = \mathrm{inv}$ and $\field = GF(4)$ (and, of course, $\alpha = \mathrm{id}$ and $\field = GF(2)$).
\end{remark}

The following result is a straightforward adaption of a result of \cite{BH/PivotLoopCompl/09}.
\begin{proposition} [Theorem~8 of \cite{BH/PivotLoopCompl/09}] \label{prop:pu_loopc}
Let $A$ be a principally unimodular $V \times V$-matrix over a field $\field$ of characteristic $2$. Then, for all $X \subseteq V$, $\mathcal{M}_{A+X} = \mathcal{M}_A+X$.
\end{proposition}
\begin{Proof}
It suffices to show the result for $X = \{j\}$ with $j \in V$. Let $Z\subseteq V$. We compare
$\det A[Z]$ with $\det (A+\{j\})[Z]$. First assume that $j\notin Z$. Then
$A[Z] = (A+\{j\})[Z]$, thus $\det A[Z] = \det (A+\{j\}) [Z]$. Now
assume that $j\in Z$, which implies that $A[Z]$ and $(A+\{j\})[Z]$ differ
in exactly one position: $(j,j)$. We may compute determinants by
Laplace expansion over the $j$-th column, and summing minors. As
$A[Z]$ and $(A+\{j\})[Z]$ differ at only the matrix-element $(j,j)$,
these expansions differ only in the inclusion of minor $\det
A[Z\setminus\{j\}]$. Thus $\det (A+\{j\})[Z] = \det A[Z] + \det
A[Z\setminus\{j\}]$, and this computation is in $GF(2)$ as $A$ is PU and $\field$ of characteristic $2$. From this the statement follows.
\end{Proof}

The following result is an adaption of \cite[Theorem~8.2]{BH/NullityLoopcGraphsDM/10} (\cite[Theorem~8.2]{BH/NullityLoopcGraphsDM/10} shows that the family of binary \dmatroids is vf-closed).
\begin{theorem} \label{thm:inv_repr_gf4_vfclosed}
The family of \dmatroids $\mathrm{inv}$-representable over $GF(4)$ is vf-closed.
\end{theorem}
\begin{Proof}
Let $M$ be a \dmatroid $\mathrm{inv}$-representable over $GF(4)$. Then $M = \mathcal{M}_A*X$ for some $\mathrm{inv}$-symmetric $V \times V$-matrix $A$ over $GF(4)$ and $X \subseteq V$. Let
$\varphi$ be a sequence of vertex flips over $V$.
Let $W \in \mathcal{M}_A*X \varphi$, and consider now $\varphi'
= *X\varphi*W$. By the $S_3^V$ group structure of vertex flips (see \cite[Theorem~12]{BH/PivotLoopCompl/09}),
$\varphi'$ can be put in the following normal form:
$\mathcal{M}_A \varphi' = \mathcal{M}_A+Z_1*Z_2+Z_3$ for some
$Z_1,Z_2,Z_3 \subseteq V$ with $Z_1 \subseteq Z_2$. By Theorem~\ref{thm:inv_symm_PU}, $A$ is PU. By Proposition~\ref{prop:pu_loopc},
$\mathcal{M}_A + Z_1 = \mathcal{M}_{A+Z_1}$. Thus
$\mathcal{M}_A+Z_1*Z_2+Z_3 = \mathcal{M}_{A+Z_1}*Z_2+Z_3$. By
construction $\emptyset \in \mathcal{M}_A \varphi'$. Hence we
have $\emptyset \in \mathcal{M}_{A+Z_1}*Z_2$. Therefore $Z_2
\in \mathcal{M}_{A+Z_1}$ and so $A+Z_1*Z_2$ is defined.
Consequently, $A' = A+Z_1*Z_2+Z_3$ is defined and
$\mathcal{M}_A \varphi' = \mathcal{M}_{A'}$. Hence $M\varphi =
\mathcal{M}_A*X\varphi = \mathcal{M}_{A'}*W$ and thus $\mathrm{inv}$-symmetric matrix $A'$
represents $M\varphi$. Consequently, $M\varphi$ a \dmatroid $\mathrm{inv}$-representable  over $GF(4)$.
\end{Proof}

In contrast with Theorem~\ref{thm:inv_repr_gf4_vfclosed}, it is shown in \cite{BH/NullityLoopcGraphsDM/10} that there are \dmatroids $\mathrm{id}$-representable over $GF(4)$ that are \emph{not} \vfsafe.

\section{Quaternary matroids and bicycle matroids}

Let $M = (V, \mathcal{B})$ be a matroid representable over $\field$, and described by its bases. Let $B$ be a standard representation of $M$ over $\field$. Then $B$ is equal to
$$
\bordermatrix{&X&V\setminus X\cr
X & I & S
}
$$
for some $X \in \mathcal{B}$, where $I$ is the identity matrix of suitable size. Let $\alpha$ be an automorphism of $\field$ that is an involution. We define $R(B,\alpha)$ to be the $\alpha$-symmetric $V \times V$-matrix
$$
\bordermatrix{&X&V\setminus X\cr
X & 0 & S\cr
V\setminus X & \alpha(-S^T) & 0
}
$$

The next result is from \cite{bouchet1987}.
\begin{proposition} [Theorem~4.4 of \cite{bouchet1987}] \label{prop:bouchet_twist_matroid}
Let $M$ be a matroid representable over $\field$, $B$ be a $X \times V$-matrix over $\field$ that is a standard representation of $M$. Then $\mathcal{M}_A = M*X$ with $A = R(B,\mathrm{id})$.
\end{proposition}

Note that if $A = R(B,\alpha)$ and $A'= R(B,\mathrm{id})$, then $A[Y]$ is nonsingular iff $A'[Y]$ is nonsingular for all $Y \subseteq V$. Hence $\mathcal{M}_A = \mathcal{M}_{A'}$.

Hence, by Proposition~\ref{prop:bouchet_twist_matroid}, if a matroid $M$ is ($\mathrm{id}$-)representable over $\field$, then $M$ is $\alpha$-representable for all automorphisms $\alpha$ of $\field$ that are involutions. Conversely, if $M$ is $\alpha$-representable for involution and automorphism $\alpha$ of $\field$, then by the exact same reasoning as the only-if direction of the proof of Theorem~4.4 of \cite{bouchet1987}, we have that matroid $M$ is representable over $\field$.
Consequently, a matroid $M$ is representable over $\field$ in the usual (matroid) sense iff $M$ is $\alpha$-representable for some involution and automorphism $\alpha$ of $\field$ iff $M$ is $\alpha$-representable for all automorphisms $\alpha$ of $\field$ that are involutions.
Therefore, choosing $\alpha = \mathrm{id}$ may not necessarily be the most natural extension of the matroid notion of representability to \dmatroids. Indeed, in view of Theorem~\ref{thm:inv_repr_gf4_vfclosed} and the remark below it, we argue that over $GF(4)$, $\mathrm{inv}$-representability is the most natural extension of the matroid notion of representability to \dmatroids.

In particular, by Proposition~\ref{prop:bouchet_twist_matroid} every quaternary matroid is a \dmatroid $\mathrm{inv}$-representable over $GF(4)$. Hence by Theorem~\ref{thm:inv_repr_gf4_vfclosed} we have the following result, which was conjectured in \cite{BH/NullityLoopcGraphsDM/10}.
\begin{corollary}
Every quaternary matroid is \vfsafe.
\end{corollary}

In this paragraph use terminology of \cite{DBLP:journals/jct/Vertigan98}. 
Let $L \subseteq \field^V$ be a subspace of $\field^V$, and let $\alpha$ be an automorphism of $\field$. We define $\mathrm{bd}(L,\field,\alpha) = \dim(L\cap\alpha(L^{\perp}))$, where $L^{\perp}$ is the orthogonal subspace of $L$.
If $|\field| = q \in \{2,3,4\}$, then there is an automorphism $\alpha: \field \rightarrow \field$ with $\alpha(x) = x^{-1}$ for all $x \in \field \setminus \{0\}$. 
In these cases, $\mathrm{bd}(L,\field,\alpha)$ is called the \emph{bicycle dimension} of $L$ and we denote it by $\mathrm{bd}(L,q)$.
Let $L$ be a subspace of $GF(4)^V$. In line with \cite{DBLP:journals/jct/Vertigan98}, we call $L \cap \mathrm{inv}(L^{\perp})$ the \emph{bicycle space} of $L$, and denote it by $\mathcal{BC}_L$. We know from \cite{DBLP:journals/jct/Vertigan98} that the dimension of $\mathcal{BC}_L$ is determined by $M(L)$. We now extend this result by showing that the \emph{matroid} of $\mathcal{BC}_L$ is determined by $M(L)$. Moreover we give an explicit formula for $M(\mathcal{BC}_L)$. Also, the proof of this result below is direct, and therefore not obtained as a consequence of an evaluation of the Tutte polynomial as in \cite{DBLP:journals/jct/Vertigan98}.

\begin{theorem} \label{thm:quaternary_bicycle_matroid}
Let $M$ be a quaternary matroid, and let $A$ be a representation of $M$ over $GF(4)$. Then the matroid $M(\mathcal{BC}_{\ker(A)})$ is equal to $\max(M+V)$.
\end{theorem}
\begin{Proof}
Let $A'$ be a standard representation of $M$ with $L = \ker(A') = \ker(A)$. Let $P = R(A',\mathrm{inv})$. By Proposition~\ref{prop:bouchet_twist_matroid}, $\mathcal{M}_P*X = M$ for some $X \subseteq V$.
By Theorem~\ref{thm:inv_symm_PU}, $P$ is PU and by Proposition~\ref{prop:pu_loopc}, $\mathcal{M}_P+V = \mathcal{M}_{P+V}$. Now, $A' = M+V[X,V]$. Recall from, e.g., \cite[Proposition~2.2.23]{Oxley/MatroidBook-2nd}, that $\ker(I \quad S)^{\perp} = \ker(-S^T \quad I)$. Thus, $\ker(P+V[V\setminus X,V]) = \mathrm{inv}(L^{\perp})$. Consequently, $\ker(P+V) = \mathcal{BC}_{L}$. Now, $M(\mathcal{BC}_{L}) = M(P+V) = \max(\mathcal{M}_{P+V}) = \max(\mathcal{M}_P+V) = \max(M*Y+V) = \max(M+V\dual Y) = \max(M+V)$.
\end{Proof}

While $\mathcal{BC}_{A}$ and $\mathcal{BC}_{A'}$ may differ when $A$ and $A'$ are different representations of $M$ over $GF(4)$, Theorem~\ref{thm:quaternary_bicycle_matroid} shows that $M(\mathcal{BC}_{A}) = M(\mathcal{BC}_{A'})$.

A binary matroid $M$ has an odd number of bases iff the dimension
of the bicycle space is zero---a result originally shown in \cite{Chen/VectorSpaceGraph/1971} for the case where $M$
is a graphic matroid. Theorem~\ref{thm:quaternary_bicycle_matroid} shows that this statement holds in general for quaternary matroids $M$.


\bibliography{../geneassembly}

\begin{thebibliography}{1}

\bibitem{bouchet1987}
A.~Bouchet.
\newblock Representability of {$\Delta$}-matroids.
\newblock In {\em Proc. 6th Hungarian Colloquium of Combinatorics, Colloquia
  Mathematica Societatis J\'{a}nos Bolyai}, volume~52, pages 167--182.
  North-Holland, 1987.

\bibitem{DBLP:conf/ipco/Bouchet95}
A.~Bouchet.
\newblock Coverings and delta-coverings.
\newblock In E.~Balas and J.~Clausen, editors, {\em IPCO}, volume 920 of {\em
  Lecture Notes in Computer Science}, pages 228--243. Springer, 1995.

\bibitem{BH/PivotLoopCompl/09}
R.~Brijder and H.J. Hoogeboom.
\newblock The group structure of pivot and loop complementation on graphs and
  set systems.
\newblock {\em European Journal of Combinatorics}, 32:1353--1367, 2011.

\bibitem{BH/NullityLoopcGraphsDM/10}
R.~Brijder and H.J. Hoogeboom.
\newblock Nullity and loop complementation for delta-matroids.
\newblock To appear in SIAM Journal on Discrete Mathematics, preprint
  [arXiv:1010.4497], 2013.

\bibitem{Chen/VectorSpaceGraph/1971}
W.-K. Chen.
\newblock On vector spaces associated with a graph.
\newblock {\em SIAM Journal on Applied Mathematics}, 20:526--529, 1971.

\bibitem{Kodiyalam_Lam_Swan_2008}
V.~Kodiyalam, T.Y. Lam, and R.G. Swan.
\newblock Determinantal ideals, {P}faffian ideals, and the principal minor
  theorem.
\newblock In {\em Noncommutative Rings, Group Rings, Diagram Algebras and Their
  Applications}, pages 35--60. American Mathematical Society, 2008.

\bibitem{Oxley/MatroidBook-2nd}
J.G. Oxley.
\newblock {\em Matroid theory, Second Edition}.
\newblock Oxford University Press, 2011.

\bibitem{tucker1960}
A.W. Tucker.
\newblock A combinatorial equivalence of matrices.
\newblock In {\em Combinatorial Analysis, Proceedings of Symposia in Applied
  Mathematics}, volume~X, pages 129--140. American Mathematical Society, 1960.

\bibitem{DBLP:journals/jct/Vertigan98}
D.~Vertigan.
\newblock Bicycle dimension and special points of the {T}utte polynomial.
\newblock {\em Journal of Combinatorial Theory, Series B}, 74(2):378--396,
  1998.

\end{thebibliography}

\end{document}